\documentclass[11pt]{amsart}

\newtheorem{theorem}{Theorem}[section]
\newtheorem{lem}[theorem]{Lemma}

\theoremstyle{definition}

\theoremstyle{remark}

\numberwithin{equation}{section}

\newcommand{\Hor}{{\mathcal{H}}}
\newcommand{\V}{{\mathcal{V}}}
\newcommand{\ra}{\rightarrow}

\newcommand{\Ad}{\text{Ad}}

\newcommand{\lb}{\langle}
\newcommand{\rb}{\rangle}

\newcommand{\mg}{\mathfrak{g}}
\newcommand{\mh}{\mathfrak{h}}
\newcommand{\mk}{\mathfrak{k}}
\newcommand{\mm}{\mathfrak{m}}
\renewcommand{\mp}{\mathfrak{p}}
\newcommand{\bs}{\backslash}

\newcommand{\R}{\mathbb{R}}

\begin{document}

\newcommand{\spacing}[1]{\renewcommand{\baselinestretch}{#1}\large\normalsize}
\spacing{1.14}

\title{Obstruction to positive curvature on homogeneous bundles}

\author {Kristopher Tapp$^\ast$}

\address{Department of Mathematics\\ Williams College\\
Williamstown, MA 01267}
\email{ktapp@williams.edu}
\thanks{$^\ast$Supported in part by NSF grant DMS--0303326.}



\begin{abstract}
Examples of almost-positively and quasi-positively curved spaces of the form $M=H\bs((G,h)\times F)$ were discovered recently~\cite{Wilking},\cite{QT}.  Here $h$ is a left-invariant metric on a compact Lie group $G$, $F$ is a compact Riemannian manifold on which the subgroup $H\subset G$ acts isometrically on the left, and $M$ is the orbit space of the diagonal left action of $H$ on $(G,h)\times F$ with the induced Riemannian submersion metric.  We prove that no new examples of strictly positive sectional curvature exist in this class of metrics.  This result generalizes the case $F=\{\text{point}\}$ proven by Geroch~\cite{Geroch}.
\end{abstract}

\maketitle


\section{Introduction}\label{intro}
In~\cite{Geroch}, Geroch studied metrics of the form $M=H\bs(G,h)$, where $h$ is a left-invariant metric on a compact Lie group $G$, $H\subset G$ is a Lie subgroup, and $M$ is the quotient (the right coset space) with the induced Riemannian submersion metric.  He proved that no new examples of positive curvature could be found among such metrics; more precisely, if $H\bs(G,h)$ has positive sectional curvature, then $H\bs G$ admits a normal homogeneous metric of positive curvature.  These ``Geroch metrics'' are generally inhomogeneous, although his problem was motivated by the well-known classification of homogeneous spaces with positive curvature.

We consider examples of the more general form $M=H\bs((G,h)\times F)$, where $F$ is a compact Riemannian manifold on which $H$ acts isometrically on the left, so $M$ is the orbit space of the diagonal left action of $H$ on $(G,h)\times F$ with the induced Riemannian submersion metric.  Topologically, $M$ is the total space of a homogeneous $F$-bundle over $H\bs G$.  If the action of $H$ on $F$ is transitive with isotropy group $K\subset H$, then $M$ is diffeomorphic to $G/K$.

Generalizing Geroch metrics in this way enriches the family of examples.  For example, the positively curved non-normal homogeneous spaces discovered in~\cite{AW} can be re-described in this form; see Section 4.  Among Geroch metrics, only the normal homogeneous examples appear.  Further, examples of this form exist with quasi-positive and almost-positive curvature~\cite{Wilking},\cite{QT}.  It is not known whether new quasi-positive curvature examples exist among Geroch metrics.

Our main result says there are no new examples of this form with positive curvature:
\begin{theorem}\label{T:main} If $M=H\bs((G,h)\times F)$ has positive sectional curvature, and $\dim(H)<\dim(G)$, then $M$ admits a homogeneous metric with positive sectional curvature.
\end{theorem}

The dimension hypothesis disallows $G=H$, in which case $M=G\bs((G,h)\times F)$ is diffeomorphic to $F$.  This method of changing the metric on $F$ via an isometric $G$-action is curvature non-decreasing; see~\cite{Cheeger}.

The author is pleased the thanks Frank Morgan for useful discussions about this work.


\section{Curvature formulas}
In this section, we derive curvature formulas for $M=H\bs((G,h)\times F)$.  Let $h_0$ denote a bi-invariant metric on $G$ (and also its restriction to $H$).  In this section, we assume that $F$ is the normal homogeneous space $F=(H,h_0)/K$ for some subgroup $K\subset H$.  This assumption is not necessary in Theorem~\ref{T:main}, but it makes our formulas cleaner.  Let $\mk\subset\mh\subset\mg$ denote the Lie algebras of $K\subset H\subset G$.

The value of $h$ at the identity $e\in G$ is determined in terms of $h_0$ by some self-adjoint $\Phi:\mg\ra\mg$ defined so that for all $A,B\in\mg$,
$$h(A,B)=h_0(\Phi(A),B).$$
Define $\mm=\mh\ominus\mk$ and $\mp=\mg\ominus\mh,$ where ``$\ominus$'' means the orthogonal compliment with respect to $h_0$.  Also denote $\lb A,B\rb := h_0(A,B)$, $|A|^2:=h_0(A,A)$ and $|A|_h:=h(A,A)$ for $A,B\in\mg$.  Let $\pi:(G,h)\times F\ra M$ denote the quotient map.  The metric on $M$ is defined so that $\pi$ is a Riemannian submersion.

\begin{lem}\label{L} Let $g\in G$ and $y\in H$.  The horizontal space of $\pi$ at $(g,yK)$ is:
$$\Hor_{(g,yK)}=\{(dL_{g}(\Phi^{-1}(\Ad_{(g^{-1}y)}u)),-dL_y(u^{\mm}))\mid u\in\mg\ominus\mk\},$$
where $u^{\mm}$ denotes the $h_0$-orthogonal projection of $u$ onto $\mm$.
\end{lem}

Above, the tangent space to $F=(H,h_0)/K$ at the coset $yK$ is identified with the horizontal space at $y$ of $(H,h_0)\ra F$, which is $dL_y(\mm)$.

\begin{proof}
The dimension of our purported horizontal space is the dimension of $\mg\ominus\mk$, which is correct because $M$ is diffeomorphic to $G/K$.  The vertical space of $\pi$ is:
$$\V_{(g,yK)}=\{(dR_{g}A,p(dR_y(A)))\mid A\in\mh\},$$
where $p$ denotes the projection onto $dL_y(\mm)$.  It remains to verify that $\Hor$ and $\V$ are orthogonal.  Letting $A\in\mh$ and $u\in\mg\ominus\mk$,
\begin{eqnarray*}
\lefteqn{\lb (dL_{g}(\Phi^{-1}(\Ad_{(g^{-1}y)}u)),-dL_y(u^{\mm})),(dR_{g}A,p(dR_y(A)))\rb_M}\hspace{1in} \\
  & = & h(\Phi^{-1}(\Ad_{(g^{-1}y)}u),\Ad_{g^{-1}}A) + \lb-dL_y(u^{\mm}),dR_y(A)\rb \\
  & = & \lb\Ad_{(g^{-1}y)}u,\Ad_{g^{-1}}A\rb - \lb\Ad_y(u^{\mm}),A\rb \\
  & = & \lb\Ad_yu,A\rb - \lb\Ad_y(u^{\mm}),A\rb = \lb\Ad_y(u-u^{\mm}),A\rb = 0
\end{eqnarray*}
\end{proof}

Any element of $M$ has the form $\pi(g_0,eK)$ for some $g_0\in G$.  For $u\in\mg\ominus\mk$, let $\overline{u}$ denote the element of $\Hor_{(g_0,eK)}$ corresponding to $u$ as in Lemma~\ref{L}:
$$\overline{u} = (dL_{g_0}(\Phi^{-1}(\Ad_{g_0^{-1}}u)),-u^{\mm})$$

For $X,Y\in\mg\ominus\mk$, let $k_M(g_0,X,Y)$ denote the unnormalized sectional curvature of $d\pi(\overline{X})$ and $d\pi(\overline{Y})$, which by ONeal's formula equals:
\begin{equation}\label{curv}k_M(g_0,X,Y) = k_h(\Phi^{-1}(\Ad_{g_0^{-1}}X),\Phi^{-1}(\Ad_{g_0^{-1}}Y)) + k_F(X^\mm,Y^\mm) + (3/4)|[\overline{X},\overline{Y}]^{\V}|^2.\end{equation}
Here $k_h$ and $k_F$ denote the unnormalized sectional curvatures of $(G,h)$ at $g_0$ and $F$ at $eK$ respectively.  The point $\pi(g_0,eK)$ has positive curvature if and only if $k_M(g_0,X,Y)>0$ for all linearly independent $X,Y\in\mg\ominus\mk$.

The $k_h$ term can be expanded by P\"uttmann's formula from~\cite{put}:
\begin{gather}\label{E:put}
k_h(Z_1,Z_2) = (1/2)\lb[\Phi(Z_1),Z_2] + [Z_1,\Phi(Z_2)],[Z_1,Z_2]\rb - (3/4)|[Z_1,Z_2]|^2_h \\
  +\lb B(Z_1,Z_2),\Phi^{-1}(B(Z_1,Z_2))\rb - \lb B(Z_1,Z_1),\Phi^{-1}(B(Z_2,Z_2))\rb,\notag
\end{gather}
where $B(Z_1,Z_2) = (1/2)([Z_1,\Phi(Z_2)]+[Z_2,\Phi(Z_1)])$.

In the Lie bracket term of Equation~\ref{curv}, $\overline{X}$ and $\overline{Y}$ are assumed to be extended to horizontal vector fields on $(G,h)\times F$ in a neighborhood of $(g_0,eK)$.  The natural extensions suggested by Lemma~\ref{L} are:
$$\overline{X}(g,yK) = (dL_{g}(\Phi^{-1}(\Ad_{(g^{-1}y)}X)),-dL_y(X^{\mm})),$$
and similarly for $\overline{Y}$.  This definition of $\overline{X}(g,yK)$ depends on the coset representative $y$, but becomes well-defined once we choose a section $\mathcal{S}$ of $(H,h_0)\ra F$ in a neighborhood of $e$ from which to choose our $y$'s.  In other words, $\mathcal{S}$ is an open submanifold of $H$ whose tangent space at $e$ is $\mm$, which projects to a neighborhood of $eK$ in $F$.  With these extensions, we have:
\begin{lem}\label{A}At the point $(g_0,eK)$,
\begin{gather*}
[\overline{X},\overline{Y}] = \big(dL_{g_0}([\Phi^{-1}(\Ad_{g_0^{-1}}X),\Phi^{-1}(\Ad_{g_0^{-1}}Y)]
                              +\Phi^{-1}(\Ad_{g_0^{-1}}([Y^\mm,X]-[X^\mm,Y])) \\
                              - \Phi^{-1}[\Phi^{-1}(\Ad_{g_o^{-1}}X),\Ad_{g_0^{-1}}Y] +
                              \Phi^{-1}[\Phi^{-1}(\Ad_{g_0^{-1}}Y),\Ad_{g_0^{-1}}X]),-[X^\mm,Y^{\mm}]^\mm\big).
\end{gather*}
\end{lem}
\begin{proof}
Let $\overline{A}$ denote the following vector field on $(G,h)\times F$:
$$\overline{A}(g,yK) = (dL_g(\Phi^{-1}(\Ad_{g_0^{-1}} X)),-p(dR_y(X^\mm))).$$
Notice that $\overline{A}$ agrees with $\overline{X}$ at $(g_0,eK)$, but $\overline{A}$ is easier to work with because it's a product of a left-invariant field on $G$ and a Killing field on $F$.  Define $\overline{B}$ to be the analogous product vector field that agrees with $\overline{Y}$ at $(g_0,eK)$.  Then at $(g_0,eK)$,
$$[\overline{X},\overline{Y}]   = [\overline{A}+(\overline{X}-\overline{A}),\overline{B}+(\overline{Y}-\overline{B})]
                                = [\overline{A},\overline{B}] + [\overline{A},(\overline{Y}-\overline{B})]
                                                            + [(\overline{X}-\overline{A}),\overline{B}].$$

For the first term,
$$[\overline{A},\overline{B}] = (dL_{g_0}[\Phi^{-1}(\Ad_{g_0^{-1}}X),\Phi^{-1}(\Ad_{g_0^{-1}}Y)],-[X^\mm,Y^\mm]^\mm).$$

The second term, $[\overline{A},(\overline{Y}-\overline{B})]$, simplifies because $\overline{Y}-\overline{B}$ vanishes at $(g_0,eK)$.  To see how, let $\{V_i\}$ be an $h_0$-orthonormal frame of left-invariant fields on $G$.  Let $\{W_j\}$ be a left-invariant extension to $\mathcal{S}$ of a $h_0$-orthonormal basis of $\mm$.  The $W_j$'s cannot generally be made everywhere tangent to $\mathcal{S}$, but they project to a local orthonormal frame on $F$ near $eK$.  Choose a path in $G\times F$ in the direction of $A(g_0,eK)$, which will have the form $t\mapsto(g_0a(t),y(t))$ where $a(t)$ is a path in $G$ with $a(0)=e$ and $a'(0)=\Phi^{-1}(\Ad_{g_0^{-1}}X)$, and $y(t)$ is a path $\mathcal{S}$ with $y'(0)=-X^\mm$.  Then,
\begin{eqnarray*}
[\overline{A},(\overline{Y}-\overline{B})]
  & = & [\overline{A},\sum_i \lb\overline{Y}-\overline{B},(V_i,0)\rb(V_i,0) + 
        \sum_j \lb\overline{Y}-\overline{B},(0,W_j)\rb(0,W_j)] \\
  & = & \sum_i\overline{A} \lb\overline{Y}-\overline{B},(V_i,0)\rb(V_i,0)
        +\sum_j\overline{A} \lb\overline{Y}-\overline{B},(0,W_j)\rb(0,W_j)\\
  & = & \sum_i\frac{d}{dt}\Big|_{t=0}\lb\Phi^{-1}(\Ad_{((g_0a(t))^{-1}y(t))}Y - \Ad_{g_0^{-1}}Y),V_i\rb V_i\\
  &  & \hspace{1in} + \sum_j\frac{d}{dt}\Big|_{t=0}\lb-Y^\mm +\Ad_{y(t)^{-1}}Y^\mm,W_j\rb W_j\\
  & = & \sum_i \lb\Phi^{-1}(\Ad_{g_0^{-1}}[y'(0),Y] + [-a'(0),\Ad_{g_0^{-1}}Y]),V_i\rb V_i \\
  &  & \hspace{1in}  + \sum_j\lb[-y'(0),Y^\mm],W_j\rb W_j
\end{eqnarray*}
which shows that:
\begin{eqnarray*}
[\overline{A},(\overline{Y}-\overline{B})]
& = & (\Phi^{-1}(\Ad_{g_0^{-1}}[y'(0),Y] + [-a'(0),\Ad_{g_0^{-1}}Y]),[-y'(0),Y^\mm]^\mm) \\
& = & (-\Phi^{-1}(\Ad_{g_0^{-1}}[X^\mm,Y] + [\Phi^{-1}(\Ad_{g_0^{-1}}X),\Ad_{g_0^{-1}}Y]),[X^\mm,Y^\mm]^\mm)
\end{eqnarray*}
Similarly,
$$[\overline{B},(\overline{X}-\overline{A})]
= (-\Phi^{-1}(\Ad_{g_0^{-1}}[Y^\mm,X] + [\Phi^{-1}(\Ad_{g_0^{-1}}Y),\Ad_{g_0^{-1}}X]),[Y^\mm,X^\mm]^\mm)$$
Collecting terms completes the proof.
\end{proof}


\section{Proof of main Theorem}
In this section, we prove Theorem~\ref{T:main} by an argument similar to Geroch's main proof in~\cite{Geroch}.  We carry over all notation from the previous section.
\begin{proof}
Suppose first that the action of $H$ on $F$ is not transitive.  For any point $p_0\in F$, there exists a $\pi$-horizontal zero-curvature plane at $(e,p_0)$; namely any plane of the form $\text{span}\{(A,0),(0,V)\}$, where $A\in\mg$ is $h$-orthogonal to $\mh$, and $V\in T_{p_0}F$ is orthogonal to the $H$-orbit.  Further, if $p_0$ is contained in a principal orbit, then the $A$-tensor of $\pi$ vanishes on this plane, since $(A,0)$ and $(0,V)$ extend to local $\pi$-horizontal fields with everywhere vanishing $F$ and $G$ components respectively.  Thus, if the action of $H$ on $F$ is not transitive, then $M$ does not have positive curvature.

Next, suppose that the action of $H$ on $F$ is transitive, so $M$ is diffeomorphic to $G/K$.  The homogeneous bundle $F\ra M\ra H\bs G$ is called ``fat'' if $[X,Y]\neq 0$ for all non-zero $X\in\mm$ and $Y\in\mp$.  We will prove:
\begin{equation}\label{con}\text{If $M$ has positive curvature, then the bundle is fat.}\end{equation}

Proving assertion~\ref{con} suffices to complete our proof of Theorem~\ref{T:main} because of Berard Bergery's classification of fat homogeneous bundles~\cite{BB},\cite{fat}.  If $(G,H)$ is a rank one symmetric pair, he proved that fatness is equivalent to the existence of a positively curved homogeneous metric on $M$.  Further, he proved that if the dimension of $F$ is greater than $1$, then fatness implies that $(G,H)$ must be a rank one symmetric pair.

So it remains to verify that if $F$ is $1$-dimensional and $M$ has positive curvature, then $M$ admits a homogeneous metric of positive curvature.  Since the projection $M\ra H\bs(G,h)$ is a Riemannian submersion, if $M$ has positive curvature, then so does $H\bs(G,h)$, which by Geroch's Theorem implies that $G/H$ admits a normal homogeneous metric of positive curvature.  It is known which circle bundles over rank one symmetric spaces admit positive curvature, and those which admit positive curvature admit homogenous metrics of positive curvature.  The three non-symmetric positively curved normal homogeneous spaces are all odd-dimensional, making $M$ even dimensional; since our metric on $M=H\bs((G,h)\times S^1)$ admits a free isometric $S^1$-action (induced by the $S^1$-action on the second factor of $(G,h)\times S^1$), $M$ could not have positive curvature because of Berger's theorem, which says that a positively curved even-dimensional manifold does not admit a nonvanishing Killing field.

Assume that the bundle is not fat, so there exists non-zero vectors $\mathcal{X}\in\mm$ and $\mathcal{Y}\in\mp$ with $[\mathcal{X},\mathcal{Y}]=0$.  We must prove that $M$ does not have positive curvature.  Since $\mathcal{X}$ and $\mathcal{Y}$ commute, they are together contained in some maximal abelian subalgebra $\mathfrak{t}\subset\mg$.  Almost every element of $\mathfrak{t}$ is generic, i.e., only commutes with other elements of $\mathfrak{t}$.  Let $X_0,Y_0\in\mathfrak{t}$ be generic element arbitrarily close to $\mathcal{X}$ and $\mathcal{Y}$.  Notice that $[X_0,Y_0]=0$.  We will freely use ``$\epsilon$'' to denote any quantity that goes to zero as $X_0\ra\mathcal{X}$ and $Y_0\ra\mathcal{Y}$.  For example, $|[Y_0^\mm,X_0] - [X_0^\mm,Y_0]|<\epsilon$,  and $|X_0^\mm,Y_0^\mm|<\epsilon$, which is significant because terms of these forms appear in Lemma~\ref{A}.

Define $f:G\ra\R$ as:
$$f(g) =\lb \Ad_{g^{-1}}X_0,\Phi^{-1}(\Ad_{g^{-1}}X_0)\rb.$$  Let $g_0\in G$ denote a global maximum of $f$.  We will prove that $k_M(g_0,X_0,Y_0)\leq\epsilon$.  Since $M$ is compact, this will establish that $M$ does not have positive curvature.  Let $X=\Ad_{g_0^{-1}}X_0$ and $Y=\Ad_{g_0^{-1}}Y_0$.  Notice that $X$ and $Y$ are generic, and $[X,Y]=0$.

Since $g_0$ is a critical point of $f$, we have for all $Z\in\mg$,
$$0=
\frac{d}{dt}\Big|_{t=0}\lb\Ad_{(g_0e^{tZ})^{-1}}X_0,\Phi^{-1}(\Ad_{(g_0e^{tZ})^{-1}}X_0)\rb
  = 2\lb[-Z,X],\Phi^{-1} X\rb = 2\lb[\Phi^{-1} X,X],Z\rb, 
$$
which shows that:
\begin{equation}\label{e1}[\Phi^{-1} X,X]=0.\end{equation}

Since the second derivative of $f$ is nonpositive at $g_0$ along any path, for all $Z\in\mg$:
\begin{eqnarray}\label{second}
0 & \geq & (1/2)\frac{d^2}{dt^2}\Big|_{t=0}\lb\Ad_{(g_0e^{tZ})^{-1}}X_0,\Phi^{-1}(\Ad_{(g_0e^{tZ})^{-1}}X_0)\rb\\
& = & \frac{d}{dt}\Big|_{t=0}\lb[-Z,\Ad_{(g_0e^{tZ})^{-1}}X_0],\Phi^{-1}(\Ad_{(g_0e^{tZ})^{-1}}X_0)\rb\notag\\
& = & \lb[-Z,X],\Phi^{-1}[-Z,X]\rb + \lb[-Z,[-Z,X]],\Phi^{-1} X\rb\notag\\
& = & \lb[Z,X],\Phi^{-1}[Z,X]\rb - \lb[Z,\Phi^{-1} X],[Z,X]\rb.\notag
\end{eqnarray}
It follows from equation~\ref{second} that if $Z$ commutes with $\Phi^{-1} X$ then $Z$ commutes with $X$.  Since $X$ is generic, the converse holds: if $Z$ commutes with $X$ then $Z$ commutes with $\Phi^{-1} X$.  In particular, 
\begin{equation}\label{e2}[Y,\Phi^{-1} X]=0.\end{equation}

Equation~\ref{curv} says that:
$$k_M(g_0,X_0,Y_0) = k_h(\Phi^{-1}X,\Phi^{-1}Y) + k_F(X_0^\mm,Y_0^\mm) + (3/4)|[\overline{X_0},\overline{Y_0}]^{\V}|^2.$$
For the second term, $k_F(X_0^\mm,Y_0^\mm)<\epsilon$.  For the first term, using Equations~\ref{E:put},\ref{e1},\ref{e2}:
\begin{eqnarray*} k_h(\Phi^{-1}X,\Phi^{-1}Y) & = &
   (1/2)\lb[X,\Phi^{-1}Y],[\Phi^{-1}X,\Phi^{-1}Y]\rb \\ 
 & & -(3/4)|[\Phi^{-1}X,\Phi^{-1}Y]|^2_h +(1/4)\lb[X,\Phi^{-1}Y],\Phi^{-1}[X,\Phi^{-1}Y]\rb.
\end{eqnarray*}
For the third term, we temporarily add the hypothesis that $F$ is normal homogeneous, as in the previous section, in which case Lemma~\ref{A} and Equations~\ref{e1} and~\ref{e2} yield:
\begin{eqnarray*}
(3/4)|[\overline{X_0},\overline{Y_0}]^{\V}|^2 & < &
  (3/4)|[\Phi^{-1}X,\Phi^{-1}Y]-\Phi^{-1}[X,\Phi^{-1}Y]|^2_h + \epsilon\\
  & = & (3/4)|[\Phi^{-1}X,\Phi^{-1}Y]|^2_h+(3/4)|\Phi^{-1}[X,\Phi^{-1}Y]|^2_h \\
  &   & -(3/2)\lb[\Phi^{-1}X,\Phi^{-1}Y],[X,\Phi^{-1}Y]\rb+\epsilon
\end{eqnarray*}
Combining terms gives:
$$k_M(g_0,X_0,Y_0)<\lb[\Phi^{-1} Y,X],\Phi^{-1}[\Phi^{-1} Y,X]\rb - \lb[\Phi^{-1} Y,\Phi^{-1} X],[\Phi^{-1} Y,X]\rb
  +\epsilon <\epsilon.$$
The final inequality is justified by substituting $Z=\Phi^{-1}Y$ into Equation~\ref{second}.

It remains to handle the case where $F$ is non-normal homogeneous.  In this case, $F$ can be expressed as $F=(H,h')/K$, where $h'$ is a left-invariant right-$K$-invariant metric on $H$.  The value of $h'$ at $e$ can be chosen to agree with $h_0$ on $\mk$ and to preserve the orthogonality of $\mk$ and $\mm$.  The value of $h'$ on $\mm$ is determined by some endomorphism $\varphi:\mm\ra\mm$ defined so that $h'(A,B)=h_0(A,\varphi B)$ for all $A,B\in\mm$.  This added generality only affects a slight change to the formulas on the previous section.  Lemma~\ref{L} changes to:
$$\Hor_{(g,yK)}=\{(dL_{g}(\Phi^{-1}(\Ad_{(g^{-1}y)}u)),-dL_y(\varphi^{-1}(u^{\mm})))\mid u\in\mg\ominus\mk\},$$
and Lemma~\ref{A} becomes:
\begin{gather*}
[\overline{X},\overline{Y}] = \big(dL_{g_0}([\Phi^{-1}(\Ad_{g_0^{-1}}X),\Phi^{-1}(\Ad_{g_0^{-1}}Y)]
                              +\Phi^{-1}(\Ad_{g_0^{-1}}([\varphi^{-1}(Y^\mm),X]-[\varphi^{-1}(X^\mm),Y])) \\
                              - \Phi^{-1}[\Phi^{-1}(\Ad_{g_o^{-1}}X),\Ad_{g_0^{-1}}Y] +
                              \Phi^{-1}[\Phi^{-1}(\Ad_{g_0^{-1}}Y),\Ad_{g_0^{-1}}X]),-[\varphi^{-1}(X^\mm),\varphi^{-1}(Y^{\mm})]^\mm\big).
\end{gather*}
The previous proof of the case $\varphi=\text{Identity}$ works equally well for arbitrary $\varphi$.
\end{proof}

\section{Metric Variations}
In this section, we consider the family
$M_t=H\bs((G,h_t)\times F)$, where $h_t$ is a family of left-invariant metric on $G$, with $h_0$ bi-invariant.  As in Section 2, we assume that $F$ is normal homogeneous.  The following lemma says that the initial metric $M_0$ is non-normal homogeneous:
\begin{lem}
$M_0=H\bs ((G,h_0)\times F)$ is isometric to $(G,\tilde{h})/K$, where $\tilde{h}$ is a left-invariant right-$H$-invariant metric on $G$.
\end{lem}
The metric $\tilde{h}$ is defined by $(G,\tilde{h}) = ((G,h_0)\times(H,h_0))/H$.  The non-normal homogeneous spaces discovered in~\cite{AW} have the form $(G,\tilde{h})/K$ for this choice of $\tilde{h}$, as described in~\cite{Eschenburg}.
\begin{proof}
First, $M_0=H\bs((G,h_0)\times((H,h_0)/K))$ is the quotient of $(G,h_0)\times(H,h_0)$ under the action of $H\times K$ defined by $(h',k')\star(g,h) = (h'g,h'hk'^{-1})$.  Second, $(G,\tilde{h})/K = ( ((G,h_0)\times (H,h_0))/H) /K$
is the quotient of $(G,h_0)\times(H,h_0)$ under the action of $H\times K$ which sends $(h',k')\star(g,h)=(gh',k'^{-1}hh')$.
Define an isometry between these two quotients of $(G,h_0)\times(H,h_0)$ by sending the orbit of $(g,h)$ to the orbit of $(g^{-1},h^{-1})$.
\end{proof}

Many examples in~\cite{Wilking} and~\cite{QT} have the form $M_t=H\bs((G,h_t)\times F)$ and have quasi- or almost-positive curvature for all $t>0$.  It is interesting that these examples are variations of homogeneous metrics.  For example, $K=SO(n-1)\subset H=SO(n)\subset G=SO(n+1)$ gives a family of almost-positively curved metrics on $T^1S^n$.  In this case, the starting non-normal homogeneous metric $M_0$ is a Levi-Civita connection metric on $T^1S^n$, which is geometrically a more natural starting point than the normal homogeneous metric.  In these examples, $(G,H)$ and $(H,K)$ are rank one symmetric pairs, and the variation is $\Phi_t(A) = (1-t)A^{\mh} + A^{\mp}$ for $A\in\mg$ and $t\in[0,1)$.  This describes a family of nonnegatively curved left-invariant metrics on $G$ obtained by steadily shrinking vectors tangent to $H$.  Differentiating the function $f(t)=k_{M_t}(g,X,Y)$ with $X,Y\in\mg\ominus\mk$ and $g\in G$ chosen so that $f(0)=0$ provides an alternative way to derive the conditions in~\cite{QT} under which points have positive curvature for $t>0$.

We omit this derivation, but point out that more than one derivative of $f$ (in fact three) are needed.  This is not surprising.  For any family $M_t$ of nonnegatively curved compact spaces with $M_0$ homogeneous, if a single point becomes positively curved to first order, then there exists a variation whereby all points become positively curved to first order, and hence the space admits strictly positive curvature.  This is because, by compactness, there is a finite collection of variations whereby every point becomes positively curved for at least one in the collection, and the first variation of curvature formula for a sum of metric variations is additive in the variations; see for example~\cite{Strake}.


\bibliographystyle{amsplain}

\end{document}